\documentclass[a4j]{article}
\newtheorem{thm}{Theorem}[section]
\newtheorem{prop}[thm]{Proposition}
\newtheorem{cor}[thm]{Corollary}

\newtheorem{defn}[thm]{Definition}
\usepackage{amssymb}

\newcommand{\e}{\varepsilon}
\newcommand{\h}{\hbar}
\newcommand{\dbar}{d\!\!\lower-0.8ex\hbox{-}}

\title{Star exponential functions as two-valued elements}
\author{
H.Omori\thanks{E-mail:omori@ma.noda.tus.ac.jp}
\\
Department of Mathematics, \\
Faculty of Science and Technology\\
Tokyo University of Science, \\
Noda, Chiba, 278-8510, Japan.
\and
Y.Maeda\thanks{E-mail:maeda@math.keio.ac.jp}
\\
Department of Mathematics, \\
Faculty of Science and Technology, \\
Keio University \\
Hiyoshi, Yokohama, 223-8522, Japan.
\and
N.Miyazaki\thanks{E-mail:miyazaki@math.hc.keio.ac.jp}
\\
Department of Mathematics, \\
Faculty of Economics,\\
Keio University, \\
Hiyoshi, Yokohama, 223-8521, Japan.
\and
A.Yoshioka\thanks{ E-mail:yoshioka@rs.kagu.tus.ac.jp}
\\
Department of Mathematics,  \\
Faculty of Science\\
Tokyo University of Science, \\
Kagurazaka,
Tokyo, 102-8601, Japan.
}

\date{}

\begin{document}
\maketitle
\begin{abstract}
We propose a relatively new notion of two-valued elements, which arises
naturally in constructing the star exponential functions of the quad-ratics
in the Weyl algebra over the complex number field.  This notion enables us
to describe the group like objects of the set of star exponential functions
of quadratics in the Weyl algebra.
\end{abstract}
\section{Introduction}
\par

Geometries are described within a framework of manifolds which
are set up among the topological spaces. The question then
may arise as to whether there are possibilities to employ
other notions rather than manifolds.
In this paper, we attempt to propose a notion of
{\it two-valued elements}, which seems to renew a geometric concept.
\par
A nontrivial example of objects we propose in this paper is given
by the Hopf-fibering $S^3\stackrel{S^1}{\rightarrow}S^2$.
Viewing $S^3=\coprod_{q\in S^2}S_q^1$ (disjoint union),
we consider the double covering $\tilde{S}_q^1$
of each fiber $\tilde{S}_q^1$, which is denoted by $\tilde{S}^3$.
When $\tilde{S}^3$ is considered as a point set, we are able
to define local trivializations of 
$\tilde{S}^3|_{V_i}\cong V_i\times \tilde{S}^1$ naturally through the
trivializations ${S}^3|_{V_i}$ given on a simple
open covering $\{V_i\}_{i\in \Gamma}$ of $S^2$. This structure permit
us to treat $\tilde{S}^3$ as a local Lie group, and hence it
looks like a topological space. On the other hand, we have
a projection
$\pi:\tilde{S}^3=\coprod_{q\in S^2}\tilde{S}_q^1\rightarrow S^3=
\coprod_{q\in S^2}{S}_q^1$
as the union of fiberwise projections, as if it were a non-trivial double
covering. $\tilde{S}^3$ cannot be a manifold, since $S^3$ is
simply connected.  This might suggest us
to make the notion of points vague. In particular, the
``points"of $\tilde{S}^3$ should be regarded as {\it two-valued elements}
with $\pm$ ambiguity.
\par
We now consider a 1-parameter subgroup $S^1$ of $S^3$ and the
inverse image $\pi^{-1}(S^1)$. Since all points of ${\widetilde S}^3$
are ``two-valued", this simply looks like a combined object of
$S^1\times {\mathbb Z}_2$ and the double covering group, i.e. in some 
restricted region, this object can be regarded as a point set by
several ways. In such a region, the ambiguity 
only arises in the case two pictures of point sets are mixed up. 
Similar phenomena appear in constructing star exponential functions
of quadratic forms in the suitably extended Weyl system, which leads us to open 
a new concept of geometry as a noncommutative (quantum) aspect.

 \par
In the paper \cite{OMMY5}, we have shown strange phenomena which break 
associativity for the Weyl algebra over ${\mathbb C}$ generated by two 
generators $u$ and $v$. Furthermore, we have shown that the Lie
algebra over $\mathbb C$ of quadratic forms can be exponentiated to 
the ``group'' which looks as if it were a double covering group 
of $SL_{\mathbb C}(2)$ which is simply connected, or the 
complexification of the metaplectic group $Mp(2,{\mathbb R})$.   

As a sequel to this work, we develop to the case of Weyl algebra with
$2m$-generators $u_1,\cdots,u_m,v_1,\cdots,v_m$, and show that similar 
phenomena occur also in the case of $2m$-generators. We show that 
star exponential functions can be viewed as two-valued elements.  
We note that an approach using the notion of gerbes will be a
possibility to describe such phenomena (cf.  \cite{Bry}, \cite{Hi}, \cite{OMMY5}),
which will also give rise to a new geometrical formulation.

\section{Weyl algebra and orderings}
\subsection{Weyl algebra}

\par

The Weyl algebra $W_\h$ is the algebra over ${\mathbb C}$ generated
by $u_1,\cdots,u_m,v_1,\cdots,v_m$ with the following commutation
relations:
\begin{equation}\label{rel}
[u_i,v_j]=-i\h\delta_{ij},
\end{equation}
where $\h$ is a positive constant and $[a,b]=a{*}b-b{*}a$.
Here, the product on $W_\h$ is denoted by $*$.
For abbreviation, we set ${\mathbf u}=(u_1,\cdots,u_m),$
and ${\mathbf z}=({\mathbf u},{\mathbf v})
=(u_1,\cdots,u_m,v_1,\cdots,v_m)$.
Let $Sym(2m,\mathbb C)$ be the set of complex symmetric matrix
$A=(A_{ij})$.
For $A\in Sym(2m,\mathbb C)$, we define a quadratic form by
\begin{equation}
A_*({\mathbf z})=\sum_{i,j=1}^{2m}A_{ij}\frac{1}{2}(z_i{*}z_j+z_j{*}z_i).
\end{equation}
Denote by ${\cal A}_\h$ the set of $A_*(z)$,
where $A\in Sym(2m,\mathbb C)$.
It is easily seen that ${\cal A}_\h$ forms a complex Lie algebra
isomorphic to$sp_{\mathbb C}(m)$.

\subsection{Orderings}
\par
Orderings are treated in the physical literature (cf. \cite{AW}) 
in quantum mechanics as the rule of association from $c$-number 
functions to $q$-number functions. There are typical orderings, called  
the standard ordering, the antistandard ordering and the Weyl
ordering, and in case of complex variables $\zeta_k=u_k{+}iv_k$, 
$\zeta_l^*=u_l{-}iv_l$, the normal ordering, 
the antinormal ordering and the Weyl ordering. 

However, from the mathematical view point, it is better to go back to 
the original understanding of Weyl which says that the ordering is 
the problem of realization of the Weyl algebra $W_\h$. 
Since the Weyl algebra is the universal enveloping algebra of 
Heisenberg Lie algebra, the Poincar{\'e}-Birkhoff-Witt theorem   
shows that this algebra can be viewed as an algebra defined on 
the space of polynomials.  

For precise formulations of ordering prescriptions in {\it formal}
deformation quantization, one can refer to the article {\cite{BNW}},
but the theory using a formal deformation parameter gives only a probe
for genuine quantum theory. We emphasize here that the deformation 
parameter $\h$ in this note is {\it not} a formal parameter, but a
parameter moving among positive reals.

Thus, we generalize orderings as follows. Let $J$ be a $2m\times 2m$ 
matrix defined by 
$J={\scriptsize{\left[\matrix{0 & -I_m \cr I_m & 0}\right]}}$. For every
symmetric complex $2m\times 2m$ matrix $K=(K^{ij})$, we set the product

\begin{equation}
\label{KK}
f({\mathbf z})*_{K}g({\mathbf z})=
f\exp\{\frac{i\h}{2}({\sum}_{i,j=1}^{2m}\overleftarrow{\partial_{z_i}}
\Gamma^{ij}\overrightarrow{\partial_{z_j}})\}g,
\end{equation}
where $\Gamma=(\Gamma^{ij})=(K^{ij}{+}J^{ij})$.
The product formula (\ref{KK}) is
well-defined for all $f,g\in {\mathbb C}[{\mathbf z}]$,
where ${\mathbb C}[{\mathbf z}]={\mathbb C}[z_1,\cdots,z_{2m}]$,
and this satisfies  
\begin{equation}\label{K-rel}
z_i{*_K}z_j-z_j{*_K}z_i(= [z_i,z_j]_{*_K})=i\h J^{ij},
\end{equation}
which give the same commutation relations (\ref{rel}) as the Weyl algebra $W_\h$.

\begin{prop}
\label{prop-realization}
For every complex symmetric $2m\times 2m$ matrix $K$,
$({\mathbb C}[{\mathbf z}],*_K)$ forms
an associative algebra isomorphic to $W_\h$.
\end{prop}

Proposition~\ref{prop-realization} gives a realization of
the Weyl algebra $W_\h$, and at the same time, it also gives the way
of the expressions of elements of the Weyl algebra $W_\h$.
For instance, computing $u_i*u_j*u_k$ by using (\ref{KK}) gives the
expression of $u_i*u_j*u_k$ as a polynomial. Thus, the product formula
(\ref{KK}) will be referred to as $K$-ordering, i.e. giving an
ordering is nothing but a product formula which gives
the Weyl algebra $W_{\h}$ where generators are fixed.
Note that according to the choice of $K$:
$$
 \left[\matrix{
   0 & I_m\cr
   I_m & 0}
 \right],\,\,
 \left[\matrix{
   0 & -I_m\cr
   -I_m & 0}
 \right],\,\,
 \left[\matrix{
   0 & 0\cr
   0 & 0}
 \right],
$$
the product formulas (\ref{KK}) gives the standard ordering,
the antistandard ordering and the Weyl ordering respectively.

By the above formulation of orderings, intertwiners between
$K$-orderings are explicitly given as follows:
\begin{prop}
For every pair of complex symmetric $2m\times 2m$ matrices $K, K'$, we have
the intertwiner 
$T_\h:({\mathbb C}[{\mathbf z}],*_K)\rightarrow({\mathbb C}[{\mathbf
z}],*_{K'})$
defined as
\begin{equation}\label{intertwiner}
T_\h(f)=\exp\Big(\frac{\h}{2i}\sum_{i,j}(K^{ij}-K^{'ij})
\partial_{z_i}\partial_{z_j}\Big)f.
\end{equation}
Namely the following identity
\begin{equation}\label{intertwiner2}
T_\h(f*_Kg)=T_\h(f)*_{K'}T_\h(g),
\end{equation}
holds for any $f,g \in {\mathbb C}[{\mathbf z}]$.
\end{prop}

Although the presice statement will be given in the forthcoming paper, 
the intertwiner can be extended for a certain class of functions.  
However, as it has be shown in \cite{OMMY5} the intertwiner 
behaves only a 2-to-2 mappings in the space of exponential functions 
of quadratic forms, since the square root appears in the amplitude 
part of intertwined functions.

We think this is a basic phenomena which breaks the associativity 
of $*$-product in the space of closed linear hull of the exponential
functions of quadratic forms. Such strange phenomena occurs only in the 
the case that the deformation parameter is a non-formal parameter. 
In spite of this, it is important that 
{\it one can consider one parameter subgroups} via the theory of
ordinary differential equations.

\section{Star exponential functions }
\subsection{Star exponential functions }
\par

We give the explicit formula for the star exponential function
$e_*^{tA_*({\mathbf z})}$ via $K$-ordering.
For $A\in Sym(2m,\mathbb C)$, we denote by $A[{\mathbf z}]$
the symmetric quadratic function defined by
\begin{equation}\label{quad-func}
A[{\mathbf z}]=\sum_{i,j=1}^{2m}A_{ij}z_i  z_j.
\end{equation}
Set ${\mathbb C}^{\times}={\mathbb C}-\{0\}$, and we denote by ${\cal F}$
the set defined by
\begin{equation}\label{cal-F}
{\cal F}=\{F=g\exp Q[{\mathbf z}]\,|\, g\in {\mathbb C}^{\times}, Q\in
Sym(2m,\mathbb C)\}.
\end{equation}
For $A\in Sym(2m,\mathbb C)$, we set as
\begin{equation}\label{circ-quadratic}
A_{*_K}({\mathbf z})=\sum_{i,j=1}^{2m}A_{ij}\frac{1}{2}(z_i{*_K}z_j+z_j{*_K}z_i).
\end{equation}
The product formula (\ref{KK}) gives
\begin{equation}\label{translate}
A_{*_K}({\mathbf z})=A[{\mathbf z}] + i\h {\rm Tr}K A.
\end{equation}
We realize the star exponential functions of $A_{*}({\mathbf z})$, for
$A\in Sym(2m,\mathbb C)$ with the help of $K$-ordering.
Namely, in order to get the formula $e_{*_K}^{tA_{*_K}({\mathbf z})}$,
we set $F_K(t)=e_{*_K}^{A_{*_K}({\mathbf z})}$, and consider the
following equation:
\begin{equation}\label{propagation}
\left\{\begin{array}{l}
\partial_tF_K(t)=A_{*_K}({\mathbf z})*_K F_K(t), \\
  \,\,\,  F_K(0)=1.
\end{array}\right.
\end{equation}
By the product formulas (\ref{KK}) and
(\ref{circ-quadratic}), the evolution equation (\ref{propagation})
can be expressed as a differential equation. Thus the uniqueness of the
real analytic solution holds, if it exists.

By setting
\begin{equation}\label{formula-solution}
F_K(t)=g_K(t)\exp Q_K(t)[{\mathbf z}],\quad
{\mbox{where}}\,\,Q_K(t)[{\mathbf z}]\in Sym(2m,\mathbb C),
\end{equation}
the evolution equation (\ref{propagation}) is reduced to a system of
ordinary differential equations on
$g_K(t)$ and $Q_K(t)[{\mathbf z}]$.
By a direct computation, although it is rather complicated, we have
\begin{thm}\label{theorem-solution}
The evolution equation (\ref{propagation})
has the unique analytic solution $F_K(t) \in {\cal F}$ explicitly given by
\begin{equation}\label{formula-solution}
F_K(t)=g_K(t)\exp Q_K(t)[{\mathbf z}],\nonumber
\end{equation}
where
\begin{equation}\label{Q-solution}
Q_K (t)=\frac{-J}{\h}(\tan \h t{J A})\cdot (I-iK\tan \h t{J A})^{-1}
\end{equation}
\begin{equation}\label{g-solution}
g_K(t)=\big(\det(\cos t\h{J A}-iK\sin t\h{J A})\big)^{-1/2}.
\end{equation}
\end{thm}
\par\medskip
\noindent{\bf Remark}\,\,
Millard \cite{m} also obtained this product formula by solving
the successive power series of a Riccati-type equation. 
Remark also that for every $t\in{\mathbb C}$ there is $K$-ordering 
such that $F_K(t)$ is well-defined.    

\par\medskip
Among the formula (\ref{formula-solution}) of the star exponential
functions of $A_{*_K}({\mathbf z})$, we particularly have for
the standard ordering and the Weyl ordering by plugging
$K=\scriptsize{\left[\begin{array}{cc} 0 & I_m\\I_m & 0
\end{array}\right]}$,
and $K=0$, respectively.
\par
In particular, we have
\begin{cor}\label{N-star-exponential-cor}
For any $A\in Sym(2m,\mathbb C)$, the star exponential function
$ e_{*}^{tA_{*}({\mathbf z})}$ is expressed as
\begin{equation}\label{N-star-exponential}
e_{*}^{tA_{*}({\mathbf z})}
=\det(\cos\h t{JA})^{-1/2}\cdot
\exp\left(\frac{-J}{\h}(\tan\h t{JA})[{\mathbf z}]\right)
\end{equation}
via the Weyl ordering.
\end{cor}

\subsection{Star exponential functions of rank one quadratics}
\par
We examine the product formula (\ref{propagation}) by
restricting the quadratics to the rank one.
\par
For ${\mathbf x}=(x_1,\cdots,x_m)$
and ${\mathbf y}=(y_1,\cdots,y_m)\in {\mathbb C}^m$,
we set $\langle{\mathbf x},{\mathbf y}\rangle
=\sum_{i=1}^mx_iy_i$.
For ${\mathbf a}, \,{\mathbf b}\in {\mathbb C}^m$, we consider
$\langle{\mathbf a},
   {\mathbf u}\rangle=\sum_{i=1}^ma_iu_i$ and
$\langle{\mathbf b},
    {\mathbf v}\rangle=\sum_{i=1}^mb_iv_i$
as elements of $W_\h$. It is easy to see
\begin{equation}\label{angle-bracket}
[\langle{\mathbf a},{\mathbf u}\rangle,
\langle{\mathbf b},{\mathbf v}\rangle]_{*}
=-i\h \langle{\mathbf a},{\mathbf b}\rangle.
\end{equation}
Hence,
if $\langle{\mathbf a},{\mathbf a}\rangle=1$, then
$\langle{\mathbf a},{\mathbf u}\rangle$ and
$\langle{\mathbf a},{\mathbf v}\rangle$ form
a canonical conjugate pair.
Let $S_{\mathbb C}^{m-1}=
\{{\mathbf a}\in {\mathbb C}^m\,|\,\langle{\mathbf a},{\mathbf
a}\rangle=1\} $.
For every ${\mathbf a}\in S_{\mathbb C}^{m-1}$, and $\alpha, \beta,
\gamma\in {\mathbb C}$,
we consider a quadratic form
\begin{equation}\label{rk-1-quadratic-form}
\begin{array}{ll}
&B_*(\alpha,\beta,\gamma)\\
&=\alpha \langle{\mathbf a},{\mathbf u}\rangle{*}
                 \langle{\mathbf a},{\mathbf u}\rangle
{+}\beta \langle{\mathbf a},{\mathbf v}\rangle{*}
                 \langle{\mathbf a},{\mathbf v}\rangle
{+}\gamma\big(\langle{\mathbf a},{\mathbf u}\rangle{*}
                 \langle{\mathbf a},{\mathbf v}\rangle
{+}\langle{\mathbf a},{\mathbf v}\rangle{*}
                 \langle{\mathbf a},{\mathbf u}\rangle\big),
\end{array}
\end{equation}
which is called a {\it quadratic form of rank one}.
\par
In the following, we assume that
the discriminant $D=\gamma^2-\alpha\beta=1$.
We now write down the star exponential for the quadratic form
$B_*(\alpha,\beta,\gamma)$ of rank one.
We denote by $F_M(\alpha,\beta,\gamma)$
and $F_N(\alpha,\beta,\gamma)$ the solution of (\ref{propagation})
for $A_*({\mathbf z})=B_*(\alpha,\beta,\gamma)$
with respect to $K=0$ and
${\scriptsize{\left[\matrix{
  0 & I_m\cr
  I_m & 0}
\right]}} $
respectively.
Then, we have (see also \cite{OMMY5}):
\begin{cor}\label{rk-1-star-exponential-solution-cor}
Assume that $D=\gamma^2-\alpha\beta=1$.
Then, for ${\mathbf a}\in {\mathbb C}^m$
such that $\langle{\mathbf a},{\mathbf a}\rangle=1$, we have
\begin{equation}\label{rk-1-M-exponential-solution-formula}
F_M(t,\alpha,\beta,\gamma)
=g_M(t,\alpha,\beta,\gamma) \cdot \exp Q_M(t,\alpha,\beta,\gamma),
\end{equation}
where
\begin{equation}\label{rk-1-M-exponential-solution-g}
g_M(t,\alpha,\beta,\gamma)=
(\cos \h t)^{-1},
\end{equation}
\begin{equation}\label{rk-1-M-exponential-solution-Q}
Q_M(t,\alpha,\beta,\gamma)
=\frac{1}{\h}(\tan \h t)\cdot
\big(\alpha\langle{\mathbf a},{\mathbf u}\rangle^2{+}
 \beta \langle{\mathbf a},{\mathbf v}\rangle^2{+}
 2\gamma\langle{\mathbf a},{\mathbf u}\rangle\langle{\mathbf a},{\mathbf
v}\rangle\big)
\end{equation}

Similarly, we have
\begin{equation}\label{rk-1-N-exponential-solution-formula}
F_N(t,\alpha,\beta,\gamma)
=g_N(t,\alpha,\beta,\gamma) \cdot \exp Q_N(t,\alpha,\beta,\gamma).
\end{equation}
Here, $g_N$ and $Q_N$ are given by
\begin{equation}\label{rk-1-N-exponential-solution-g}
g_N(t,\alpha,\beta,\gamma)
= e^{-i\h t\gamma}
\cdot\left({\cos 2\h t-i\gamma\sin 2\h t}\right)^{-1/2} ,
\end{equation}
\begin{equation}\label{rk-1-N-exponential-solution-Q}
Q_N(t,\alpha,\beta,\gamma)
=\frac{1}{\hbar}\big(X_N(t) \langle{\mathbf a},{\mathbf u}\rangle^2
+Y_N(t) \langle{\mathbf a},{\mathbf v}\rangle^2
+2Z_N(t) \langle{\mathbf a},{\mathbf u}\rangle{\circ}
\langle{\mathbf a},{\mathbf v}\rangle\big) ,
\end{equation}
where
\begin{equation}\label{N-XYZ}
\left\{\begin{array}{l}
X_N(t)=\frac{\alpha}{2}
           \left(\frac{\sin 2\h t}{\cos 2\h t -i\gamma\sin 2\h t}\right),
\\
Y_N(t)=\frac{\beta}{2}
           \left(\frac{\sin 2\h t}{\cos 2\h t -i\gamma\sin 2\h t}\right),
\\
Z_N(t)=\frac{i}{2}
           \left(1-\frac{1}{\cos 2\h t -i\gamma \sin 2\h t}\right)
\end{array}
\right.
\end{equation}
and the $\circ$ in the product simply means that we use
the standard ordering.
\end{cor}

\section{Polar elements are two-valued elements}
\subsection{Polar elements}
\par
Using the formulas of the star exponential functions
(\ref{rk-1-M-exponential-solution-formula}) and
(\ref{rk-1-N-exponential-solution-formula}),
we show how two-valued elements appear.
\par
We give justifications of the star exponential function
of quadratic forms $B_*(\alpha,\beta,\gamma)$ as follows.
We consider $B_*(\alpha,\beta,\gamma)$ defined by
(\ref{rk-1-quadratic-form}), and consider
the star exponential functions expressed by the standard ordering and the Weyl ordering.
Looking at the formulas in Corollary
\ref{rk-1-star-exponential-solution-cor}
and evaluating for $t=\frac{\pi}{2\h}$, we have that
$F_M(\frac{\pi}{2\h},\alpha,\beta,\gamma)$ diverges, however
$F_N(\frac{\pi}{2\h},\alpha,\beta,\gamma)$
has a meaning.
\par
Thus, we think of $F_{N}(\frac{\pi}{2\h},\alpha, \beta, \gamma)$ as
a realization of the star exponential function of
$\frac{\pi}{2\h} B_{*}(\alpha, \beta, \gamma)$, which is denoted by
$\exp_{*}(\frac{\pi}{2\h}B_{*}(\alpha, \beta, \gamma))$.
However, by the formula (\ref{rk-1-N-exponential-solution-formula})
in Corollary \ref{rk-1-star-exponential-solution-cor}, we obtain
\begin{thm}\label{theorem1}
Assume ${\mathbf a } \in
S_{\mathbb C}^{m-1}$. For any $(\alpha, \beta, \gamma)$
with ${\gamma}^{2}-\alpha \beta=1$, we have
\begin{equation}\label{B-exp-formula}
\exp_*\left( \frac{\pi}{2\h}B_*(\alpha,\beta,\gamma) \right)
=\sqrt{-1}e_{\circ}^{\frac{2i}{\h}\langle{\mathbf a},{\mathbf
u}\rangle{\circ}
 \langle{\mathbf a},{\mathbf v}\rangle }
\end{equation}
which is independent of the choice of $\alpha, \beta, \gamma$.
\end{thm}

We will show that the ambiguity of $\sqrt{-1}$ can not be eliminated
for all $(\alpha, \beta, \gamma)$.

\begin{defn}
Assume ${\mathbf a} \in S_{\mathbb C}^{m-1}$.
\begin{equation}\label{polar-element}
\varepsilon_{00}({\mathbf a})=\exp_*\left(\frac{\pi}{2\h}B_*(0,0,1)\right)
\end{equation}
is called the polar element.
\end{defn}
%
%
\subsection{Two-valued elements}
\par
We explain that the polar elements $\varepsilon_{00}({\mathbf a})$,
${\mathbf a}\in S_{\mathbb{C}}^{m-1}$
play the same role like the two-valued elements as below.
Since $(\alpha,\beta,\gamma)=(0,0,1)$ and $(0,0,-1)$ are arcwise connected
in the set $\gamma^2\!-\!\alpha\beta=1$,
and thus, they have to be viewed as a single element.
By Theorem~\ref{theorem1},
we have
\begin{equation}
  \label{eq:contra}
\begin{array}{lll}
e_*^{\frac{\pi}{2\h}
(\langle{\mathbf a},
   {\mathbf v}\rangle{*}
   \langle{\mathbf a},
   {\mathbf u}\rangle+ \langle{\mathbf a},
   {\mathbf u}\rangle{*}
   \langle{\mathbf a},
   {\mathbf v}\rangle )}&=&
\sqrt{-1}e_{\circ}^{\frac{2i}{\h}\langle{\mathbf a},
   {\mathbf u}\rangle{\circ}
   \langle{\mathbf a},
   {\mathbf v}\rangle}\\
&=&
   e_*^{-\frac{\pi}{2\h}(\langle{\mathbf a},
   {\mathbf v}\rangle{*}
   \langle{\mathbf a},
   {\mathbf u}\rangle+ \langle{\mathbf a},
   {\mathbf u}\rangle{*}
   \langle{\mathbf a},
   {\mathbf v}\rangle
   )}.
\end{array}
\end{equation}
However, considering the exponential law of the
$*$-exponential function
$$e_*^{\frac{t}{2\h}(
\langle{\mathbf a},
   {\mathbf v}\rangle{*}
   \langle{\mathbf a},
   {\mathbf u}\rangle+ \langle{\mathbf a},
   {\mathbf u}\rangle{*}
   \langle{\mathbf a},
   {\mathbf v}\rangle )}$$
for
$t\in {\mathbb C}-\{$singular set$\}$, we must set
\begin{equation}
\begin{array}{lll}
e_*^{\frac{\pi}{2\h}(
\langle{\mathbf a},
   {\mathbf v}\rangle{*}
   \langle{\mathbf a},
   {\mathbf u}\rangle+ \langle{\mathbf a},
   {\mathbf u}\rangle{*}
   \langle{\mathbf a},
   {\mathbf v}\rangle )}
&=&ie_{\circ}^{\frac{2i}{\h}
\langle{\mathbf a},
   {\mathbf u}\rangle{\circ}
   \langle{\mathbf a},
   {\mathbf v}\rangle
},\\
e_*^{-\frac{\pi}{2\h}(
\langle{\mathbf a},
   {\mathbf v}\rangle{*}
   \langle{\mathbf a},
   {\mathbf u}\rangle+ \langle{\mathbf a},
   {\mathbf u}\rangle{*}
   \langle{\mathbf a},
   {\mathbf v}\rangle
)}
&=&-ie_{\circ}^{\frac{2i}{\h}
\langle{\mathbf a},
   {\mathbf u}\rangle{\circ}
   \langle{\mathbf a},
   {\mathbf v}\rangle
}.
\end{array}
\end{equation}
If one wants to fix the sign ambiguity, the exponential law and
(\ref{eq:contra}) gives
\begin{equation}\begin{array}{l}
-1
=e_*^{\frac{\pi}{2\h}(
\langle{\mathbf a},
   {\mathbf v}\rangle{*}
   \langle{\mathbf a},
   {\mathbf u}\rangle+ \langle{\mathbf a},
   {\mathbf u}\rangle{*}
   \langle{\mathbf a},
   {\mathbf v}\rangle
)}
*e_*^{\frac{\pi}{2\h}(
\langle{\mathbf a},
   {\mathbf v}\rangle{*}
   \langle{\mathbf a},
   {\mathbf u}\rangle+ \langle{\mathbf a},
   {\mathbf u}\rangle{*}
   \langle{\mathbf a},
   {\mathbf v}\rangle
)}\\
\quad=e_*^{\frac{\pi}{2\h}(
\langle{\mathbf a},
   {\mathbf v}\rangle{*}
   \langle{\mathbf a},
   {\mathbf u}\rangle+ \langle{\mathbf a},
   {\mathbf u}\rangle{*}
   \langle{\mathbf a},
   {\mathbf v}\rangle
)}*
e_*^{-\frac{\pi}{2\h}(
\langle{\mathbf a},
   {\mathbf v}\rangle{*}
   \langle{\mathbf a},
   {\mathbf u}\rangle+ \langle{\mathbf a},
   {\mathbf u}\rangle{*}
   \langle{\mathbf a},
   {\mathbf v}\rangle
)}
=1.
\end{array}
\end{equation}
\par
We choose a continuous path of
$(\alpha,\beta,\gamma)$ from $(0,0,1)$ to $(0,0,-1)$ for
the case $m=1$ concretely as follows: Set
$\langle{\mathbf a},{\mathbf u}\rangle=u,
\langle{\mathbf a}, {\mathbf v}\rangle=v$ and ${\e}_{00}$
stands for
${\e}_{00}(\mathbf a)$.
By a careful calculation, we see
\par
\noindent
\setlength{\unitlength}{.5mm}
\begin{picture}(80,50)(0,10)
\thicklines
\qbezier[200](50,15)(10,15)(10,50)
\qbezier[200](10,50)(10,85)(50,85)
\qbezier[200](50,85)(60,85)(60,50)
\qbezier[200](60,50)(60,15)(50,15)
\put(56,75){$\bullet$}
\put(60,77){$\e_{00}$}
\put(42,20){$\bullet$}
\put(28,20){$-\e_{00}$}
\thinlines
\put(62,50){\vector(0,-1){6}}
\qbezier[200](50,85)(40,85)(40,50)
\qbezier[200](40,50)(40,15)(50,15)
\qbezier[200](12,52)(15,75)(58,77)
\qbezier[50](12,52)(22,64)(58,65)
\qbezier[50](12,52)(20,40)(58,35)
\put(32,38){$\circledcirc$}
\put(60,56){\parbox{.13\linewidth}{\tiny{This line is identically
$\e_{00}$}}}
\put(11,50){$\bullet$}
\put(6,50){$1$}
\put(15,58){$\scriptstyle{2\theta}$}
\end{picture}
\hfill\parbox[b]{.6\linewidth}
{\begin{equation}\label{adjoint}
\begin{array}{l}\quad
\mbox{Ad}(e_*^{\frac{i\theta}{2\hbar}(u^2+v^2)})e_*^{2tuv}\\=
e_*^{t(\sin 2\theta\,\,(u^2-v^2)+\cos 2\theta \,\,2uv)}.
\end{array}\end{equation}
Since the discriminant of the quadratic form of the right hand side is
identically $1$, the right hand side is identically $\e_{00}$ for
$t=\frac{\pi}{2\hbar}$. In particular,
Ad$(e_*^{\frac{\pi i}{4\hbar}(u^2+v^2)})e_*^{\frac{\pi}{\h}uv}=\e_{00}$.
On the other hand, consider, for each $\theta$, the one parameter subgroup
Ad$(e_*^{\frac{i\theta}{2\hbar}(u^2+v^2)})e_*^{2tuv}$
with respect to $t$, $t\in [0,\frac{\pi}{2\hbar}]$.
}
\par
\noindent
We easily see that Ad$(e_*^{\frac{\pi i}{4\hbar}
(u^2+v^2)})e_*^{2tuv}=e_{*}^{-2tuv}$.
In particular,

$$
{\rm{Ad}}(e_*^{\frac{\pi i}{4\hbar}(u^2+v^2)})\e_{00}=-\e_{00}
$$
by the exponential law. Move $2\theta$ from $0$ to ${\pi}$.
Then, we see the desired fact. Note also that by
(\ref{adjoint}), (\ref{rk-1-N-exponential-solution-g}),
 there is a singularity
at $2\theta=\frac{\pi}{2}$, $t=\frac{\pi}{4\hbar}$.

We also have in the standard ordering for $D=1$,
$$
e_*^{\frac{\pi}{2\hbar}(\alpha\langle{\mathbf a},{\mathbf u}\rangle^2
{+}\beta\langle{\mathbf a},{\mathbf v}\rangle^2{+}
\gamma(\langle{\mathbf a},{\mathbf u}\rangle{*}\langle{\mathbf a},{\mathbf
v}\rangle{+}
\langle{\mathbf a},{\mathbf v}\rangle{*}\langle{\mathbf a},{\mathbf
u}\rangle)}
=\sqrt{-1}e_{\circ}^{\frac{2i}{\hbar}\langle{\mathbf a},{\mathbf u}\rangle
      \circ\langle{\mathbf a},{\mathbf v}\rangle}.
$$
By the exponential law, we see that  $\e_{00}({\mathbf a})$ satisfies
\begin{equation}
  \label{eq:e00e00}
\e_{00}({\mathbf a})^2=(\e_{00}^{-1}({\mathbf a}))^2=-1,\quad
\e_{00}({\mathbf a})*\e_{00}^{-1}({\mathbf a})=1.
\end{equation}

Therefore, we must conclude that the
sign ambiguity cannot be eliminated.
One has to set
$\e_{00}({\mathbf a})=\sqrt{-1}e_{\circ}^{\frac{2i}{\h}
\langle{\mathbf a},
   {\mathbf u}\rangle{\circ}
   \langle{\mathbf a},
   {\mathbf v}\rangle
}$ with the sign ambiguity. Similar phenomena have been discussed 
by Olver \cite{Ol}.   

\par\bigskip\noindent
{\it By the above observation, the polar element $\e_{00}({\mathbf a})$ should
be regarded as a two-valued element. Otherwise we do have a contradiction $1=-1$.}
\par\bigskip\noindent

Only this way one can permit the identity $-\e_{00}({\mathbf
a})=\e_{00}({\mathbf a})$.
But since such a notion does not exist in the set theory,
it is impossible to define ${\e}_{00}({\mathbf a})$ as a point in a point
set.

\bigskip

In what follows, we set
$$
\e_{00}(k) = e_*^{\frac{\pi}{2\hbar}(u_k{*}v_k{+}v_k{*}u_k)}
=\sqrt{-1}e_{\circ}^{\frac{2i}{\hbar}u_k\circ v_k}, \quad
k=1, 2, \dots, m.
$$
These are all regarded as two-valued elements. Although it is natural
to think $\e_{00}(k){*}\e_{00}(l)=\e_{00}(l){*}\e_{00}(k)$, we also have the
equality
$$
\e_{00}(k){*}\e_{00}(l)=-\e_{00}(l){*}\e_{00}(k), \quad (k\not=l)
$$
at the same time.
Hence we have $\e_{00}(k){*}\e_{00}(k)={\pm 1}$, but we see
$\e_{00}(k)^2=-1$. This is just the same as
$\{\pm 1\}\{\pm 1\}=\{\pm 1\}$, but $\{\pm 1\}^2=1$.

Hence $\e_{00}(k)^2$ behaves like an ordinary number in the extended
Weyl algebra. In spite of this,
$\e_{00}(k)$ does not behave like an ordinary number
$i$, since it is easy to see with the bumping identity (cf. \cite{OMMY5}) that
$$
{\rm{Ad}}(\e_{00}(k))u_k =-u_k,\quad {\rm{Ad}}(\e_{00}(k))v_k =-v_k.
$$
Using this we easily have
\begin{equation}
  \label{eq:henna}
{\rm{Ad}}(\e_{00}(1))(\sum^m_{i=1}b_iu_i)=-b_1u_1+\sum^m_{i=2}b_iu_i.
\end{equation}

Since every $\langle{\mathbf a},{\mathbf u}\rangle,
{\mathbf a}\in S_{\mathbb C}^{m-1}$ is translated to $u^1$ by a symplectic
transformation, we have in general the reflection w.r.t. $\mathbf a$:
\begin{equation}
  \label{eq:orikaesi}
\begin{array}{l}
{\rm{Ad}}(\e_{00}({\mathbf a}))\langle{\mathbf b},{\mathbf u}\rangle
=\langle{\mathbf b}
-2\langle{\mathbf a},{\mathbf b}\rangle{\mathbf a}, {\mathbf u}\rangle\\
{\rm{Ad}}(\e_{00}({\mathbf a}))\langle{\mathbf b},{\mathbf v}\rangle
=\langle{\mathbf b}-2\langle{\mathbf a},{\mathbf b}\rangle{\mathbf
a},{\mathbf v}\rangle
\end{array}
\end{equation}

We introduce a notion called a blurred double covering group, which is
a group like object formed by 2-valued elements (\cite{{OMMY5}}).

\medskip
\begin{thm}
${\rm{Ad}}(\e_{00}({\mathbf a}){*}\e_{00}({\mathbf b}))$ generate
$SO(m,\mathbb C)$, hence
$\{\e_{00}(\mathbf a){*}\e_{00}(\mathbf b)\}$ generate a blurred
double covering group of $SO(m,\mathbb C)$.
However, this blurred double cover has a point set picture of
$SO(m,\mathbb C)\times {\mathbb C}^{\times}$.
If ${\mathbf a}, {\mathbf b}$ are restricted in real vectors, then
${\rm{Ad}}(\e_{00}(\mathbf a){*}\e_{00}(\mathbf b))$ generate $SO(m)$, hence
$\{\e_{00}(\mathbf a){*}\e_{00}(\mathbf b)\}$ generate a
blurred double covering group of $SO(m)$, which may be viewed as $Spin(m)$.
\end{thm}

 \thanks
  {$\!\!\!\!\!{}^\S$ Partially supported by Grant-in-Aid for
  Scientific Research (\#14540092), Ministry of Education, Culture, Sports,
Science and
  Technology, Japan.\\
  \indent ${}^\dagger$ Partially supported by Grant-in-Aid for
   Scientific Research (\#15204005), Ministry of Education, Culture,
   Sports, Science and Technology, Japan.\\
 \indent ${}^{\diamond}$ Partially supported by Grant-in-Aid for
   Scientific Research (\#15740045), Ministry of Education, Culture,
   Sports, Science and Technology, Japan.\\
  \indent ${}^*$ Partially supported by Grant-in-Aid for
   Scientific Research (\#13640088),
   Ministry of Education, Culture, Sports, Science and Technology,
   Japan.}

\end{document}